\documentclass[11pt]{amsart}
\usepackage{amsmath,amssymb,amsbsy,amscd}
\usepackage[all]{xy}
\usepackage{amsfonts}
\usepackage{dsfont}
\usepackage{pxfonts}
\usepackage{color}
\usepackage[colorlinks]{hyperref}

\newtheorem{theorem}{Theorem}[section]
\newtheorem{remark}[theorem]{Remark}
\newtheorem{proposition}[theorem]{Proposition}
\newtheorem{lemma}[theorem]{Lemma}
\newtheorem{definition}[theorem]{Definition}

\newenvironment{prof}[1][proof]{\textbf{#1:} }{\ \rule{0.5em}{0.5em}}

\begin{document}
\title{On the $Z$-decomposition of some euclidian Lie algebras}

\thanks{\small{The authors specially thank Professor Mohamed Boucetta of the University of
Cadi-Ayyad of Morocco,
 for the  suggestion of this topic and the $Z$-decomposition technic to accomplish this
work.}}

 \maketitle
\begin{center}
\author{  \textbf{R. P.  Nimpa}$^{1}$,\quad \textbf{M. B. Djiadeu}$^{2}$,\quad \textbf{J. Wouafo Kamga}$^{3}$\\
\small{e-mail: $\textbf{1}.$ nimpapefoukeu@yahoo.fr,\quad
 $\textbf{2}.$ djiadeu@yahoo.fr,\quad  $\textbf{3}.$ wouafoka@yahoo.fr\\
 University of Yaounde 1, Faculty of Science, Department of
Mathematics,  P.O. Box 812, Yaounde, Republic of Cameroon.}}
\end{center}

\begin{abstract}
In this paper  we use the  $Z-$decomposition as a tool
 to find locally symmetric left invariant
Riemannian metrics on some Lie groups. For this purpose, we need to
compute the spectrum of the curvature operator. Since the study of
this spectrum is very
 difficult, we impose restriction  on the class of
Riemannian Lie groups and  their dimension. We investigate  the
$4-$dimensional connected Riemannian Lie groups whose associated Lie
algebras are $\mathbb{A}_2\oplus
2\mathbb{A}_1,\,\mathbb{A}_{3,1}\oplus \mathbb{A}_1,\,
\mathbb{A}_{3,3}\oplus \mathbb{A}_1$ and the subclasses of $3$ and
$4-$dimensional Riemannian Lie groups which are
$\mathcal{C}-$spaces.

\end{abstract}

\section{Introduction and Main results}
Locally symmetric  Riemannian metric is defined as Riemannian metric
for which the curvature tensor is invariant under  parallel
translations. It is  well known that Riemannian metric with constant
sectional curvature, product of locally symmetric Riemannian metrics
and the left invariant Riemannian metric induced by the opposite of
the Killing form on a connected compact semi-simple Lie group, are
examples of  locally symmetric Riemannian metrics. A manifold
endowed with a locally symmetric Riemannian metric is a Riemannian
locally symmetric space. Elie Cartan  in
\cite{{Ecartan1},{Ecartan2}} gave a complete classification of
Riemannian locally symmetric spaces, but he didn't give more
information about such
 metrics.

In this paper we find
 locally symmetric left invariant
Riemannian metrics on some Lie groups.

  A Lie group
$G$ together with a left invariant Riemannian metric $g$ is called a
\emph{Riemannian Lie group.} The left invariant Riemannian metric
$g$ on $G$ induces an inner product on the Lie algebra
$\mathfrak{g}$ of $G$ and conversely, any inner product on
$\mathfrak{g}$ gives rise to a unique left invariant metric on $G$.
Let $\nabla,$  $R$\, and $\langle\, , \,\rangle$ denoted the
Levi-Civita connexion , the Riemann curvature tensor associated to
$g$ \, and the inner product induced by $g$ on the Lie algebra
$\mathfrak{g}$ of $G$.
\begin{definition}\label{deflocsym}
 A left invariant Riemannian metric $g$ on a  Lie group $G$ is
locally symmetric if $\nabla R =0.$
\end{definition}
This is equivalent to saying that for \,$x,y,z,w \in \mathfrak{g}$,
\begin{equation}\label{nablacourbure3}
 \nabla_w(R(x,y)z)= R(
  \nabla_wx,y)z + R(x,\nabla_wy)z+R(x,y)\nabla_wz,
\end{equation}
where  $R(x,y)=\nabla_{[x,\,y ]}-\nabla_x\nabla_y+\nabla_y\nabla_x$.

We find what conditions on the structure constants of the Lie
algebra $\mathfrak{g}$ are needed for a left invariant Riemannian
metric $g$ on $G$ to be locally symmetric. For Lie groups of
dimension $
>3$, it is very difficult to apply Definition \ref{deflocsym}.

 Using
the $Z-$decomposition as a tool, we find some locally symmetric
metrics on some $4$-dimensional Riemannian Lie groups. In fact, the
left invariant Riemannian metric splits as a direct product of left
invariant metrics, when the Lie algebra admits a $Z-$decomposition.
Now, we formulate the following which is our first main result:
\begin{theorem}\label{mainthm1}
Let $(G,g)$ be a connected Riemannian Lie group with the associated
euclidian Lie algebra $(\mathfrak{g}, \langle\,,\,\rangle)$.
\begin{enumerate}
\item If $\mathfrak{g} =\mathbb{A}_2\oplus 2\mathbb{A}_1$, then
$g$ is locally symmetric if and only the non-zero structure constant
in an $\langle\,,\,\rangle-$orthonormal basis is
$C_{12}^{1}=a,\,a>0;$
\item If $\mathfrak{g}= \mathbb{A}_{3,1}\oplus \mathbb{A}_1$, then
$g$ is not locally symmetric;
\item if  $\mathfrak{g}= \mathbb{A}_{3,3}\oplus \mathbb{A}_1$, then $g$
is locally symmetric if the non-zero structure constants in an
$\langle\,,\,\rangle-$orthonormal base is $C_{13}^{1}=C_{23}^{2}=a,
\quad a>0.$
\end{enumerate}
\end{theorem}
We also have:
\begin{theorem}\label{mainthm2}
Let $(G,g)$ be a connected Riemannian Lie group with an harmonic
Weyl tensor, $(\mathfrak{g}, \langle\,,\,\rangle)$ the associated
euclidian Lie algebra. If $\mathfrak{g}$
  is\,one of the euclidian Lie algebras  $2\mathbb{A}_2$,\,
$\mathbb{A}_{4,6}^{\alpha, \beta}$ with $\alpha \neq 1$,\,
$\mathbb{A}_{3,9}\oplus \mathbb{A}_1$,\,
$\mathbb{A}_{3,7}^{\alpha}\oplus \mathbb{A}_1$ with $\alpha >0$ and
$\mathbb{A}_{4,12}$, then $g$ is locally symmetric.
\end{theorem}

This paper is organized as follows: section $2$ is devoted to some
basic knowledge on the curvature tensor and locally symmetric
metrics on $3-$dimension Riemannian Lie groups. In section $3$, we
recall the algorithm  of the $Z-$decomposition. In section $4$, we
prove the theorem on the $Z-$decomposition of the euclidian  Lie
algebra of a $3-$dimension Riemannian Lie groups and we prove the
main Theorems \ref{mainthm1} and \ref{mainthm2} above.

\section{Preliminaries}
Let $(G,g)$ be an $n$-dimension connected Riemannian Lie group and
denote by $ Hol$,\quad $\mathcal{H}^{inf}$  and \,$\mathcal{H}$
 the holonomy group ,
the infinitesimal holonomy group and the primitive holonomy group at
the identity element $e$ of $G$, respectively .  Their Lie algebras
are denoted by $\mathfrak{hol}$,\quad$\mathfrak{h}^{inf}$ and
$\mathfrak{h}$, respectively. Since a Riemannian Lie group is an
analytic manifold, $\mathcal{H}^{inf}\,=\,Hol$.  For more details
about the holonomy group and his subgroups, see \cite{kobayashi}.

\subsection{The curvature tensor linear operator}
Let $x,y,v,w \in \mathfrak{g}$. We  recall that on the set
$\wedge^{2}\mathfrak{g}$ of bivectors of $\mathfrak{g}$, the inner
product denoted $\langle \,,\,\rangle_{\wedge^{2}\mathfrak{g}}$ is
defined by:
\begin{equation}\label{innerbivecto}
\begin{array}{ccc}
   \langle x\wedge y, v \wedge w\rangle_{\wedge^{2}\mathfrak{g}} & = &
   \langle x,v \rangle \langle y,w \rangle- \langle x,w \rangle \langle y,v\rangle.
               \end{array}
\end{equation}
Let $\overline{R}$ be the self-adjoint operator also called the
curvature operator,  associated to the symmetric bilinear form (the
Riemannian curvature) on the euclidian space
$(\wedge^{2}\mathfrak{g},\langle
\,,\,\rangle_{\wedge^{2}\mathfrak{g}}$),
\begin{equation}\label{curoperator}
\langle\overline{R}(x\wedge y),v \wedge
    w\rangle_{\wedge^{2}\mathfrak{g}}= R(x\wedge y, v \wedge w)=\langle R(x,y)v,w\rangle=R(x,y,v,w)
\end{equation}
\begin{remark}
 If $(e_i)_{i=1 \cdots n}$ is an orthonormal basis for $\mathfrak{g}$, the
inner product $\langle \, , \,\rangle_{\wedge^{2}\mathfrak{g}}$  is
such that the bivector $(e_i \wedge e_j)_{i< j}$ is an orthonormal
basis of $\wedge^{2}\mathfrak{g}$.
\end{remark}

The relation (\ref{innerbivecto}) induces a linear isomorphism from
the set of bivectors to the set of skew-symmetric endomorphisms,
such that for $u\wedge w \in \wedge^{2}\mathfrak{g} $ and  $x\in
\mathfrak{g}$,
\begin{equation}\label{bivectorskewsym}
u\wedge w(x)=\langle u,x\rangle w-\langle w,x\rangle u.
\end{equation}

\begin{remark}\label{eigenveciden}
Using this interpretation, we identify  the image $Im
\overline{R}$\, of $\overline{R}$ to the linear space of formal
linear combination of the elements of $\left\{R(u,v),\, u,v \in
\mathfrak{g}\right\}$, denoted $S$, see \cite{szabolocal}.
\end{remark}
 For convenience, in the rest of this paper, the skew-symmetric
endomorphism associated to
 a bivector $h$ will be denoted $\widetilde{h}$.
\begin{definition}\cite{szabolocal}
Let $h$ be any arbitrary bivector.
\begin{enumerate}
\item  The decomposition
 \begin{equation}\label{jordandecom}
    \mathfrak{g}=U_0\oplus U_1\oplus \cdots \oplus U_q
\end{equation}
where $U_0$ is the kernel of $\widetilde{h}$  and  $U_k \,\, ,1\leq
k \leq q$ \, are $2-$dimensional
 invariant  real planes of the real operator
$\widetilde{h}$, is the Jordan decomposition of the Lie algebra
$\mathfrak{g}$ with respect to $h$.
\item The number $q$ in the Jordan decomposition is the \emph{rank} of $h$.
\end{enumerate}
\end{definition}

\begin{remark}
$q=\sum m_k$ where $m_k$ is the multiplicity of the non null
eigenvalue $\lambda_k=i\nu_k,\,\,\nu_k >0$ of $\widetilde{h}$.
\end{remark}
\begin{definition}\cite{szabolocal}
Let $h$ be any arbitrary bivector. The decomposition
\begin{equation}\label{darbouxdecom}
h=\sum\limits_{k=1}^{q}v_k \wedge w_k
\end{equation}
where $v_k,\,w_k \in U_k \,$ and $U_k\,$ are the  $2-$dimensional
subspaces of the Jordan decomposition, is the Darboux decomposition
or the Darboux normal form of the bivector $h$.
\end{definition}
\begin{remark}\cite{szabolocal} \label{remdarbouxdecom}
\begin{enumerate}
\item If the multiplicity $m_k$ of non-null eigenvalue
$\lambda_k=i\nu_k$ of $\widetilde{h}$ is greater than one, then the
subspace  $ H_k=U_{k_1}\oplus U_{k_2}\oplus \cdots \oplus
U_{k_{m_k}}$
 of the Jordan
decomposition is not unique. In other words  the Darboux
decomposition
  is not unique.
  \item If the eigenvalue  $\lambda_k=i\nu_k \neq 0$ is simple, then the
  bivector $v_k \wedge w_k$ corresponding to $\lambda_k$ is unique
  and
  always occurs in every Darboux decomposition of $h$.
  \end{enumerate}
\end{remark}

\begin{definition}\cite{szabolocal}\label{defeirreduciblebivec}
The eigenvector $h \in \wedge^2\mathfrak{g}$ of the curvature
operator $\overline{R}$ is said to be irreducible if  any Darboux
normal form of $h$ does not split into two nontrivial summands such
that they are also eigenvectors of $\overline{R}$.
\end{definition}

\begin{remark}\label{remireducible}
If \,$\dim \mathfrak{g}= 4$\, and \,$dim \ker \widetilde{h}=2$,\,
then $h$ is irreducible. In fact the nonnull eigenvalue $i\nu_k,
\,\nu_k>0$ of $\widetilde{h}$ is simple.
\end{remark}

\subsection{Riemannian Lie groups with harmonic Weyl tensor  or  $\mathcal{C}-$spaces}
The Ricci tensor $r$ and the scalar curvature $s$ are defined by
$r(x,y)= tr(u \longmapsto R(x,u)y)$ and $s=tr(r),$ respectively.
Dividing the tensor $R$ by the metric $g$ in the sense of
Kulkarni-Nomizu product \cite{albesse}, we obtain the Weyl tensor
and the $1-$dimensional curvature tensor $A$ (also called the
Schouten tensor) , such that
\begin{equation}\label{curvaturesrelations1}
\begin{array}{ccc}
    R= W+A \odot g,
 \quad &\qquad \qquad
A=\dfrac{1}{n-2}\left(r-\dfrac{s\,g}{2(n-1)}\right)&\qquad
\text{and}
\end{array}
\end{equation}
\begin{equation*}
A \odot
g(x,y,u,v)=A(x,u)g(y,v)+A(y,v)g(x,u)-A(x,v)g(y,u)-A(y,u)g(x,v)
\end{equation*}
for  $x,y,u,v\in \mathfrak{g}$
  For more details, see
\cite{{albesse}}

\begin{definition}\cite{{albesse},{gladunoconf}}
A Riemannian Lie group $(G,g)$ is said to be conformally flat if its
Weyl tensor is trivial.
\end{definition}

 For a conformally
flat Riemannian Lie group, the spectrum of curvature operator
$\overline{R}:\,\,\wedge^{2}\mathfrak{g}\,\,\longrightarrow\,\,\wedge^{2}\mathfrak{g}$\,\,
 is given by O.P. Gladunova  in \cite{gladunoconf}. The following
 result holds:
\begin{proposition}\cite{gladunoconf}\label{gladutheodiago}
Let $(G,g)$ be a conformally flat $n-$dimensional Riemannian Lie
group \,\, $i.e$\, $W=0$,\, and let $(e_1,e_2,\cdots,e_n)$ be an
orthonormal basis of $\mathfrak{g}$ in which the Ricci operator $r$
and the one dimensional curvature $A$ are simultaneouly
diagonalizable. Then in the basis $(e_i\wedge e_j)_{i<j},$ the
curvature operator $\overline{R}$ is diagonalizable and the spectrum
of $\overline{R}$ is $\{K_{ij},\quad i<j\}$, where
$K_{ij}=K_{\sigma}(e_i\wedge e_j)$ is the sectional curvature in the
direction  $(e_i \wedge e_j)$.
\end{proposition}
\begin{definition}\cite{{dsvoronovdivweyl},{gladunodivweyl}}
A Riemannian Lie group $(G,g)$ of dimension $\geq 4$ is called a
$\mathcal{C}$-space or a space with harmonic Weyl tensor if $div
W=0$.
\end{definition}

An euclidian $4-$dimensional Lie algebra that admits a
$Z-$decomposition, splits into a direct sum (as linear spaces), of
non trivial Lie algebras both of  dimension $2$ or  dimension $1$
and $3$. The structure constants of the Lie algebra of a
$3-$dimensional locally symmetric Riemannian Lie group satisfies the
following conditions:
\begin{remark}\label{metricdim2}
 Riemannian Lie groups of dimension $1$ and $2$ are locally locally symmetric Riemannian
spaces.
\end{remark}
\begin{lemma}\label{lemma1}
Let\, $G$\, be a connected $3-$dimensional real unimodular Lie group
with left-invariant Riemannian metric.  $(G,g)$ is a locally
symmetric Riemannian Lie group if and only if in the Lie algebra
$\mathfrak{g}$ of $G$, there exists an
$\langle\,,\,\rangle-$orthonormal basis in which  structure
constants of the Lie algebra are presented in  table $1$:
\begin{center}
\begin{tabular}{llc}
  \hline
  Lie algebra & structure constants  & restrictions \\\hline
  $\mathbb{R}^3$ &  $C_{i,j}^{k}=0$  &  \\
  $\mathbb{R}^2\rtimes \mathfrak{so}(2)$ &$ C_{1,2}^{3}=C_{3,1}^{2}=a,$ & $a\,>\,0 $\\
  $ \mathfrak{su}(2)$ & $C_{1,2}^{3}=C_{2,3}^{1}=C_{3,1}^{2}=a$ & $a\,>\,0$ \\
  \hline
\end{tabular}
\begin{center}
Table $1$
\end{center}
\end{center}
\end{lemma}
\begin{prof}
 Using orthonormal  Milnor  basis (see \cite{Jmilnor}) for unimodular euclidian Lie algebras such that the
structure constants are $C_{1,2}^{3}=a,\quad C_{2,3}^{1}=c,\quad
C_{3,1}^{3}=b $,
 the non-null components of the
Riemannian curvature tensor $R$  are:
\begin{center}
\begin{align*}
  R(e_1,e_2)e_1&=\dfrac{-2a(a-b-c)-(a-b+c)(a+b-c)}{4}\,e_2;\\
  R(e_1,e_2)e_2&=\dfrac{2a(a-b-c)+(a-b+c)(a+b-c)}{4}\,e_1;\\
 R(e_1,e_3)e_1&=\dfrac{2b(a-b+c)+(a-b-c)(a+b-c)}{4}\,e_3;\\
  R(e_1,e_3)e_3&=\dfrac{-2b(a-b+c)-(a-b-c)(a+b-c)}{4}\,e_1;\\
   R(e_2,e_3)e_2&=\dfrac{2c(a+b-c)+(a-b+c)(a-b-c)}{4}\,e_3;\\
R(e_2,e_3)e_3&=\dfrac{-2c(a+b-c)-(a-b+c)(a-b-c)}{4}\,e_2.
\end{align*}
\end{center}
For the vanishing component
$R(e_1,e_2)e_3,\,\,R(e_1,e_3)e_2,\,\,R(e_2,e_3)e_1$, if the metric
is locally symmetric, then  the direct computation of the local
symmetry condition, see equation (\ref{nablacourbure3}) in
Definition \ref{deflocsym}, yields the system:
\begin{equation}\label{unimodularsystem1}
    \begin{array}{ccc}
             (a-b)(a+b-c)^2 & = & 0 \\
             (c-a)(a-b+c)^2 & = & 0  \\
             (c-b)(a-b-c)^2 & = & 0
           \end{array}
     .
\end{equation}
$(a,b,c)$ is a solution of (\ref{unimodularsystem1}) if and only if
$(a,b,c)\in \{(0,b,b),(a,a,0),(a,0,a),(a,a,a), a ,b \in
\mathbb{R}^{>0}\}$. At least, one of the structure constants $a,b,c
\in \mathbb{R}$ is negative.
\begin{enumerate}
\item If $(a,b,c)\in \{(0,b,b),(a,a,0),(a,0,a); a ,b \in
\mathbb{R}\}$, then the curvature tensor vanish
\,$\,\,i.e\,\,R(u,v)w\,=\,0$ forall $u,v,w \in \mathfrak{g}$.
Therefore $\nabla R \,=\,0$  and the metric is locally symmetric.
\item If $(a,b,c)\,=\,(a,a,a)$ with $a\neq 0$, then
the non vanishing components of curvature tensor we have:
\begin{equation*}
    \begin{array}{ccc}
      R(e_1,e_2)e_1=\dfrac{1}{2}a^2 e_2, & R(e_1,e_3)e_2=\dfrac{1}{2}a^2 e_3 ,& R(e_2,e_3)e_2=\dfrac{1}{2}a^2 e_3, \\
      R(e_1,e_2)e_2=-\dfrac{1}{2}a^2 e_1, & R(e_1,e_3)e_3=-\dfrac{1}{2}a^2 e_1 ,& R(e_2,e_3)e_1=-\dfrac{1}{2}a^2
      e_2.
    \end{array}
\end{equation*}
By direct computation, the equality \[
    \nabla_{e_m}(R({e_i},{e_j}){e_k})=R({e_i},{e_j})\nabla_{e_m}{e_k}
    +R(\nabla_{e_m}{e_i},{e_j}){e_k}+R({e_i},\nabla_{e_m}{e_j}){e_k}
\]
 holds for $i,j,k,m \in \{1,2,3\}.$
\end{enumerate}
Therefore the metric is locally symmetric.
\end{prof}
\begin{lemma}\label{lemma2}
Let\, $G$\, be a connected $3-$dimensional real nonunimodular Lie
group with left-invariant Riemannian metric.  $(G,g)$ is a locally
symmetric Riemannian Lie group if and only if in the Lie algebra
$\mathfrak{g}$ of\, $G$, there exist an
$\langle\,,\,\rangle-$orthonormal basis in which the structure
constants of the Lie algebra are presented in table $2$:
\begin{center}
\begin{tabular}{lll}
  \hline
  Lie algebra & structure constants  & restrictions \\ \hline
  $\mathbb{R}^4$ & commutative algebra: $C_{i,j}^{k}=0$  &  \\
  $\mathfrak{g}_I$ &$C_{1,2}^{2}=C_{1,3}^{3}=a $ & $a\,>\,0 $\\
  $\mathfrak{g}_D$ &
  $C_{1,2}^{2}=C_{1,3}^{3}=a,\quad C_{1,2}^{3}=-C_{1,3}^{2}=b \quad
   \text{or}\quad C_{1,3}^{2}=a$ & $a\,>\,0,b\,>\,0$ \\
  \hline
\end{tabular}
\begin{center}
Table $2$
\end{center}
\end{center}
\end{lemma}
\begin{prof}
Using an orthonormal Milnor basis for nonunimodular euclidian Lie
algebras such that the structure constants are $C_{1,2}^{2}=a,\quad
C_{1,2}^{3}=b,\quad C_{1,3}^{2}=c,\quad C_{1,3}^{3}=d$   with $a+d
\neq 0$ and $ac+bd = 0$.
 The non null  components of the curvature  tensor $R$ are:
\begin{equation*}
\begin{array}{cc}
  R(e_1,e_2)e_1\,=\,-(a^2+\dfrac{3}{4}b^2-\dfrac{1}{4}c^2+\dfrac{1}{2}b\,c)e_2
  ,&
  R(e_1,e_2)e_2\,=\,(a^2+\dfrac{3}{4}b^2-\dfrac{1}{4}c^2+\dfrac{1}{2}b\,c)e_1,\\
 R(e_1,e_3)e_1\,=\,-(d^2-\dfrac{1}{4}b^2+\dfrac{3}{4}c^2+\dfrac{1}{2}b\,c)\,e_3,&
  R(e_1,e_3)e_3\,=\,(d^2-\dfrac{1}{4}b^2+\dfrac{3}{4}c^2+\dfrac{1}{2}b\,c)\,e_1\\
   R(e_2,e_3)e_2\,=\,(\dfrac{1}{4}(b+c)^2-a\,d)\,e_3,&
R(e_2,e_3)e_3\,=\,-(\dfrac{1}{4}(b+c)^2-a\,d)\,e_2.
\end{array}.
\end{equation*}
if the metric is locally symmetric, then  the direct computation of
the local symmetry condition of equation (\ref{nablacourbure3}) in
Definition \ref{deflocsym}, yields the system:

\begin{equation}\label{nonunisystem2}
    \begin{array}{ccc}
             (b-c)(a^2+b^2-c^2-d^2) & = & 0 \\
             (b+c)(a^2+b^2-ad+bc) & = & 0  \\
             d(a^2+b^2-ad+bc)^2 & = & 0\\
             a(c^2+d^2-ad+cb)&=&0\\
             (b+c)(c^2+d^2-ad+bc)&=&0\\
             ac+bd&=&0\\
                 a+d&\neq&0
           \end{array}.
\end{equation}
Using computer system Maple, the set of non trivial and real
solutions of the system (\ref{nonunisystem2}) is

 $\{(a,b,-b,a),(0,0,0,d),(a,0,0,0); a , d \in
\mathbb{R}^{\ast},\,\,b \in \mathbb{R}\} $. The structure constants
$a,\,b,\,c,\,d$ are uniquely determined, if we normalize by
requiring that $a\geq d,\,b\geq c$ and $a+d >0$. Therefore the set
of non trivial solution of system (\ref{nonunisystem2}) is
$\left\{(a,b,-b,a),(a,0,0,0), a>0, b>0 \right\}$
 \end{prof}

\section{The $Z-$decomposition of the Lie algebra of Riemannian Lie groups}
\begin{definition}\cite{szabolocal}
The $V-$decomposition of the Lie algebra $\mathfrak{g}$  of a
connected Riemannian Lie group $(G,g)$ is an orthogonal and
irreducible decomposition of $\mathfrak{g}$ with respect to the
primitive holonomy group at $e\in G$.
\end{definition}

Let $h$ be an irreducible eigenvector of $\overline{R}$ with non
null eigenvalue.   $\widetilde{h}$ leaves the invariant and
irreducible subspaces  of the orthogonal $V-$decomposition
 invariant. Thus, we have the following proposition:
\begin{proposition}\cite{szabolocal}
 Let $h$ be an irreducible eigenvector of $\overline{R}$ with
 non-null eigenvalue, the non trivial invariant
 subspace
  \begin{equation*}
    H_k=U_{k_1}\oplus U_{k_2}\oplus \cdots \oplus U_{k_{m_k}}
\end{equation*}
of $\widetilde{h}$ is contained in a single invariant and
irreducible subspace of the $V-$decomposition.
\end{proposition}
\begin{remark}
\begin{enumerate}
\item The dimension of each invariant and irreducible subspace is at
leat $2$;
\item We can choose a decomposition (\ref{jordandecom}) in such a way that any
subspace $U_i, \,i>0$, is contained in  a single invariant and
irreducible subspace of the $V-$decomposition, \,\, see
\cite{szabolocal};
\item One can choose a complete system of linearly independent irreducible
eigenvectors of \,$\overline{R}$\, which form a basis of\,
$\wedge^{2}\mathfrak{g}$.
\end{enumerate}
\end{remark}
 Let $\{h_1,h_2,\cdots,h_{\rho},h_{\rho+1}, \cdots
,h_{\frac{n(n-1)}{2}}\}$ be such a system and assume that just the
first\, $\rho$\,\,\,vectors are corresponding to non-null
eigenvalues. Let $h \in S$, $ h=
\sum\limits_{i=1}^{\rho}a_i\lambda_i \widetilde{h_i},\,\,\,a_i \in
\mathbb{R} $
 and $\lambda_i$ the eigenvalue of the eigenvector $h_i$. Therefore,
$S=\text{Span}\{\widetilde{h}_1,\widetilde{h}_2, \cdots
    ,\widetilde{h}_{\rho}\}$, and $\mathfrak{h}$ is the free Lie
algebra on $\{\widetilde{h}_1,\widetilde{h}_2, \cdots
    ,\widetilde{h}_{\rho}\}$.

 For irreducible bivectors $h_k$,\, $1\leq k\leq \rho$, let us
consider the Jordan decomposition
\begin{equation*}
    \mathfrak{g}=U_{k0}\oplus U_{k1} \oplus U_{k2} \oplus \cdots
    \oplus U_{kN_k}
\end{equation*}
where $U_{k0}$ is the kernel of $\widetilde{h}_k$,\, $U_{kl} ,\,l
\in \{1,2, \cdots ,N_k\} \subset \mathbb{N}$ the  real
$\overline{h}_k$ invariant $2-$plan in the Jordan decomposition and
$N_k$ the rank of $h_k$. Let \,$H_{k}^{0}=U_{k0}$ and
$H_{k}^{1}=U_{k1} \oplus U_{k2} \oplus \cdots
    \oplus U_{kN_k}$,
\begin{equation*}
    \widetilde{h}_kH_{k}^{0}=0,\,\,\qquad \,\text{and} \,\,\,\,\qquad\widetilde{h}_kH_{k}^{1} \subset
    H_{k}^{1}.
\end{equation*}
We now construct the irreducible $V-$decomposition
\subsection{The main construction}\label{vdecompconstruc} Let us choose an arbitrary
vector $h_{k_{1}}$,\, $1\leq k_1\leq \rho$,\, and  consider its
subspaces\, $H_{k_{1}}^{0}$\, and \, $H_{k_{1}}^{1}$\,constructed
above.
\begin{enumerate}
\item If for any\, $h_i$\,, $i\neq k_1$\, the relation \, $H_{k_{1}}^{1} \subseteq
H_{i}^{0}$\, holds, then:
\begin{itemize}
\item[i)] $\mathcal{H}H_{k_{1}}^{1}\subseteq
H_{k_{1}}^{1}$ so that $H_{k_{1}}^{1}$ is a nontrivial
$\mathcal{H}-$invariant subspace of $V_{j};$
 \item[ii)] $\mathcal{H}$ acts irreducibly  on $V_{j}$.
\end{itemize}
 Therefore
 $H_{k_{1}}^{1}\,=\,V_{j},\,j>0$  is one of the subspaces of the $V-$decomposition.
\item If there exist vectors $h_{k_{2}},h_{k_{3}}, \cdots ,
h_{k_{l}}$ \, with \,$1\leq k_i \leq \rho$,\, such that
\,$H_{k_{1}}^{1} \nsubseteq H_{k_{i}}^{0}$\, holds for each $k_i$,
\begin{enumerate}
\item either\, $H_{k_{1}}^{1}+H_{k_{2}}^{1}+ \cdots
+H_{k_{l}}^{1}$\, is one of the subspace  $V_j,\,j>0$\, of the
$V-$decomposition, \label{maxsysteeigen} if for any \, $h_i$\,, \,$i
\notin \{ k_2,k_3, \cdots ,k_l\}$, the relation $
H_{k_{1}}^{1}+H_{k_{2}}^{1}+ \cdots +H_{k_{l}}^{1} \subseteq
H_{i}^{0}\,$ holds.
\item or, there exist  $h_i,$\quad $i
\notin \{ k_2,k_3, \cdots ,k_l\}$ such that $
H_{k_{1}}^{1}+H_{k_{2}}^{1}+ \cdots +H_{k_{l}}^{1} \nsubseteq
H_{i}^{0}\,$. Therefore, we extend the  system $
H_{k_{1}}^{1}+H_{k_{2}}^{1}+ \cdots +H_{k_{l}}^{1}\,$ with the
elements \,$H_{i}^{1}$.
\end{enumerate}
 This process can be repeated for  $
H_{k_{1}}^{1}+H_{k_{2}}^{1}+ \cdots +H_{k_{l}}^{1}+H_{i}^{1}$ until
we obtain the \emph{maximal system} $h_{k_{1}},h_{k_{2}},
\cdots,h_{k_{j}}$ such that for every index $i, \,1\leq i\leq j$,
there is another index $i',\, 1\leq i'< i$\, with\,
$H_{k_{i'}}^{1}\nsubseteq H_{k_{i}}^{0}$. Furthermore for every
index $i \notin \{k_1,k_2, \cdots ,k_j\}$ we get the relation
$H_{k_{s}}^{1} \subseteq H_{i}^{0} $ \, $s \in \{1,2, \cdots ,
j\}$.\quad
   $ V_{k}\,=\,H_{k_{1}}^{1}+H_{k_{2}}^{1}+ \cdots
    +H_{k_{j}}^{1},\,\,\,k>0$

is one of the invariant subspaces of the $V-$decomposition.
\item By continuing  the procedure we can construct all other
invariant subspaces $V_{j},\,j>0$.
\end{enumerate}

 The following formulas where the notation $\nabla
_{V_i}V_j \subset \,V_k$ means that for any $u_i \in V_i$\,
$\nabla_{u_i}u_j \in V_k$, holds for the spaces $V_i$:
\begin{equation}\label{covariantderivativeVi}
    \begin{array}{ccc}
      \nabla
_{V_0}V_0 \subset \,V_0, & \nabla _{V_0}V_i \subset \,V_i & \nabla
_{V_i}V_i \subset \,V_0\,+\,V_i \\
      \nabla
_{V_i}V_0 \subset \,V_0\,+\,V_i & \nabla _{V_i}V_j \subset \,V_j
&\,\,\texttt{if} \,\,i\neq j,\,\,i;\,j\neq 0.
    \end{array}
\end{equation}
Therofore, $V_0$ and $V_0+V_i$ are involutive,  see
\cite{szabolocal}.

\begin{remark}\label{vdecopremark}
\begin{enumerate}
\item If $0$ is the only eigenvalue of the linear curvature tensor
$\overline{R}$, then the curvature tensor vanishes.  The primitive
holonomy group is trivial
$i.e\,\,\,\mathcal{H}_p=\{id_{\mathfrak{g}}\}$. Therefore,
$\mathfrak{g}$ does not admit a $V-$decomposition.
\item Let $h$ be an irreducible eigenvector   associated to nonnull eigenvalue
  $\lambda$ of $\overline{R}$. If $\widetilde{h}$ is a linear automorphism of $\mathfrak{g}$, then:
\begin{itemize}
\item[i)] The dimension  of \,$\mathfrak{g}$\, is even;
\item[ii)]  The subspaces involved in the\, $V-$decomposition
are\, $H^{0}\,=\,\ker \widetilde{h}=\,\{0\}$
\,and\,\,$H^1\,=\,U_1\oplus \cdots \cdots \oplus U_q=\mathfrak{g}$.
Since $h$ is irreducible, $\mathfrak{g}$ is contained in a single
invariant subspace $V_j,\,j>0$. Therefore $\mathfrak{g}$ does not
admit a $V-$decomposition.
\end{itemize}
\end{enumerate}
\end{remark}
\subsection{The Z-decomposition}
Let $V_i,\,i > 0$  be a linear subspace  given by the $V-$
decomposition . Let us consider for $i>0$, the subspaces
\begin{equation*}
Z_{i}=\text{Span}\{ u_{1},\nabla_{u_1}u_2,\nabla_{u_1}\nabla
    _{u_2}u_3, \cdots, \nabla_{u_1}\nabla
    _{u_2} \cdots \nabla_{u_l}u_{l+1}, u_i \in V_i, l\in \mathbb{N}\setminus \{0\}\}
\end{equation*}
of $\mathfrak{g}$ and $Z_{0}$ the complete subspace in
$\mathfrak{g}$ which is totally orthogonal to the space
$Z_{1}+Z_{2}\,+ \cdots +\,Z_{k}.$

\begin{remark}
$\begin{array}{ccc}
  Z_{0} \subset V_{0},\,\, &\,\, V_{i}\,\subset\,Z_{i}&\,\,\, \text{and}\,\,\,
  \,\,\mathcal{H}^{inf}Z_{i}\subset Z_{i}.
\end{array}$
\end{remark}

 We have the
orthogonal splitting
\begin{equation}\label{zdecom}
\mathfrak{g}=Z_{0} \oplus Z_{1}\oplus Z_{2}\,\oplus \cdots \oplus
    \,Z_{k}.
\end{equation}
 see \cite{szabolocal}. It is a decomposition of the Lie algebra
    $\mathfrak{g}$ into invariant and irreducible subspaces with respect to the
     infinitesimal holonomy group at $e$.
\begin{definition}\cite{szabolocal}

The splitting (\ref{zdecom}) is called a $Z-$decomposition of the
Lie algebra $\mathfrak{g}$ of a connected Riemannian Lie group
$(G,g)$.
\end{definition}

\begin{proposition}\cite{szabolocal}
The subspaces $Z_i$ induce on $G$ a totally parallel distribution.
Thus they are involutive, and the integral manifolds are totally
geodesic.
\end{proposition}
\begin{remark}
\begin{enumerate}
\item If $(G,g)$ is locally symmetric  Riemannian Lie group,
then $V_i=Z_i$, since $\mathcal{H}^{inf}=\mathcal{H}$ and $V_i$ is a
nonnull  $\mathcal{H}-$invariant subspace of
$\mathcal{H}-$irreducible subspace $Z_i$;
\item   $Z_j,\,\,j\geq 0$ is a
Lie subalgebra of $\mathfrak{g}$;
\item If the Lie algebra $\mathfrak{g}$ admits a $Z-$decomposition, then the connected
Riemannian Lie group $(G,g)$ splits into a product of Riemannian Lie
groups, see \cite{antonioscala}.
\end{enumerate}
\end{remark}

\section{Application of $Z-$decomposition.}
We apply the $Z-$decomposition technics on the Lie algebras of
$3-$dimensional connected Lie groups, to the euclidian Lie algebras
$\mathbb{A}_2\oplus 2\mathbb{A}_1,\,\mathbb{A}_{3,1}\oplus
\mathbb{A}_1$ ,\, $\mathbb{A}_{3,3}\oplus \mathbb{A}_1$, to
conformally flat $4-$dimensional Riemannian Lie groups and to
harmonic Weyl tensor $4-$dimensional Riemannian Lie groups.

\subsection{The euclidian $3-$dimensional Lie algebras}
\begin{lemma}\label{3dinmzdecom}
Let $(G,g)$ be a $3-$dimensional Riemannian Lie group. If at least
$2$ eigenvalues of the linear curvature tensor operator are nonnull,
then the Lie algebra $\mathfrak{g}$ of $G$ does not admit a
$Z-$decomposition.
\end{lemma}
\begin{prof}
For a $3-$dimensional Riemannian Lie group, any orthonormal Milnor
basis $(e_1,e_2,e_3)$  diagonalizes the Ricci operator, see
\cite{Jmilnor} and the Weyl tensor vanishes ,see \cite{albesse}.
Therefore the matrix of the curvature operator in the basis $(e_1
\wedge e_2,e_1 \wedge e_3,e_2 \wedge e_3)$ is the diagonal matrix
$\text{diag}\{K_{12},K_{13},K_{23}\}$, where $K_{ij}$ with $i<j$, is
the sectional curvature in the direction $(e_i \wedge e_j)$, see
Proposition \ref{gladutheodiago}. The eigenvectors $h_1=e_1 \wedge
e_2,\, h_2=e_1 \wedge e_3$\, and \,$h_3=e_2 \wedge e_3$\, associated
 to the eigenvalues $K_{12},K_{13}$ and $K_{23},$ respectively, are
irreducible.

 Suppose that $\mathfrak{g}$ admits a $Z-$decomposition.
Then $\mathfrak{g}=V_0\oplus V_1$ with $\dim V_0=1$\,and \, $\dim
V_1=2$.  Let $h_k \,\text{and}\,\,h_l,\,\,k\neq l$ be the
eigenvector associated  to $2$ nonnull eigenvalues of
$[\overline{R}]$. Then $H_{k}^{1}=H^{1}_{l}$ since $H_{s}^{1}\subset
V_1,\,\,s \in\{k,l\} $ \,and $\dim H_{s}^{1}=\dim V_1$. Therefore
$h_k \,\,\text{and} \,\,h_l$ are linearly dependent and
\,$\mathfrak{g}$\, does not admit a $Z-$decomposition.
\end{prof}
\begin{theorem}
Let $(G,g)$ be a $3-$dimensional  Riemannian Lie group. The Lie
algebra $\mathfrak{g}$ of $G$ admits a $Z-$decomposition if and only
there exist on $\mathfrak{g}$ an orthonormal basis with respect to
$\langle\,, \,\rangle$ in which the nonnull structure constant is \,
$C_{12}^{2}=a,\quad a>0$. The metric $g$ is locally symmetric.
\end{theorem}
\begin{prof}
\begin{enumerate}
 \item \textbf{Unimodular $3-$dimensional euclidian Lie algebra}:

  Let
$\mathbb{B}=(e_1,e_2,e_3)$ be an orthonormal Milnor basis,
$C_{12}^{3}=a,\,C_{23}^{1}=c\,\,\text{and}\,\, C_{31}^{2}=b$ the
structure constants with respect to $\mathbb{B}$. The matrix of the
linear curvature operator in the basis $(e_1 \wedge e_2,e_1 \wedge
e_3,e_2 \wedge e_3)$ is the diagonal  matrix
$\text{diag}\{K_{12},K_{13},K_{23}\}$ where
\begin{align*}
K_{12}&=\frac{-2a(a-b-c)-(a-b+c)(a+b-c)}{4}\\
 K_{13}&=\frac{2b(a-b+c)+(a-b-c)(a+b-c)}{4}\\
 K_{23}&=\frac{2c(a+b-c)+(a-b+c)(a-b-c)}{4}
\end{align*}
are the sectional curvature of the metric $g$.
\begin{enumerate}
\item If $K_{ij}\neq 0,\,i<j$, then $\mathfrak{g}$ does not admit a
$Z-$decomposition, see Lemma \ref{3dinmzdecom}.
\item If two of the eigenvalues of $\overline{R}$ vanish,\,\,
$i.e \,\,(a,b,c)\in \{(a,a,0),(a,0,a),(0,b,b), \,a,b \in
\mathbb{R}\}$, \,\,then
 the third one also vanishes and $\mathfrak{g}$ does not admit a $Z-$decomposition.
\item If one of the eigenvalues of $\overline{R}$ vanishes, and the two other do not,
 $ i.e\,\,(a,b,c)\notin \{(a,a,0),(a,0,a),(0,b,b), \,a,b \in
\mathbb{R}\}$, then \,$\mathfrak{g}$\, does not admit a
$Z-$decomposition, see Lemma \ref{3dinmzdecom}.
\end{enumerate}

\item \textbf{Nonunimodular $3-$dimensional euclidian Lie algebra}:

 Let
$\mathbb{B}=(e_1,e_2,e_3)$ be an orthonormal Milnor basis and let
$C_{12}^{2}=a,\,\,C_{12}^{3}=b,\,\, C_{13}^{2}=c,
\,\text{and}\,\,C_{13}^{3}=d$,\,with\, $a+d > 0,\quad
ac+bd=0$,\quad$a\geq d,\quad \text{and}\quad b\geq c$ be the
normalized structure constants with respect to $\mathbb{B}$. The
matrix of the curvature operator in the basis $(e_1 \wedge e_2,e_1
\wedge e_3,e_2 \wedge e_3)$ is the diagonal matrix
$\text{diag}\{K_{12},K_{13},K_{23}\}$ where
\begin{align*}
K_{12}&=-(a^2+\frac{3}{4}b^2-\frac{1}{4}c^2+\frac{1}{2}bc)<0, \quad \text{see the scale change in \cite{Jmilnor}};\\
 K_{13}&=-(d^2-\frac{1}{4}b^2+\frac{3}{4}c^2+\frac{1}{2}bc)\\
 K_{23}&=\frac{1}{4}(b+c)^2-ad
\end{align*}
are the sectional curvature of the metric $g$.
\begin{enumerate}
\item  If $(K_{13},K_{23})\neq (0,0)$, then, by  Lemma \ref{3dinmzdecom}, the Lie algebra $\mathfrak{g}$
doest not admit a $Z-$decomposition.
\item If  $K_{13}=K_{23}=0,$ then the unique solution
    that fulfills the normalizing conditions is $(a,b,c,d)=(a,0,0,0)$.
\end{enumerate}
\end{enumerate}
The matrix of the curvature operator is
 $\text{diag}\{-a^2,0,0\}$. The
$Z-$decomposition is
 $\mathfrak{g}=Z_0\oplus Z_1$
  where
$Z_0=\text{Span}\{e_3\}$ and $Z_1=\text{Span}\{e_1,e_2\}$. The
nonnull constant structure is $C_{12}^{2}=a, \quad a>0.$
\end{prof}
\begin{remark}
Let $\mathfrak{g}=Z_0\oplus Z_1$ be a $Z-decomposition$ of a
$3-$dimensional nonunimodular euclidian Lie algebra , then
$[Z_0,Z_1]\subset Z_1.$
\end{remark}

\subsection{Proof of theorem \ref{mainthm1}}

\begin{lemma}\cite{klepikov}\label{lemklepikov}
For an arbitrary inner product $\langle\,,\rangle\,$ on each of the
Lie algebras  $\mathbb{A}_2\oplus
2\mathbb{A}_1,$\,$\mathbb{A}_{3,1}\oplus
\mathbb{A}_1,\mathbb{A}_{3,3}\oplus \mathbb{A}_1$, there exist an
$\langle\,,\rangle-$orthonormal basis with non-zero structure
constants given by the following table.
\begin{center}
\begin{tabular}{llc}
  \hline
  Lie Algebra & Structure onstants & Restriction \\ \hline
   $\mathbb{A}_2\oplus 2\mathbb{A}_1$& $C_{1,2}^{1}\,=a,\quad C_{1,2}^{4}=b$ & $a>0$ \\ \\
   $\mathbb{A}_{3,1}\oplus
\mathbb{A}_1$& $C_{2,3}^{1} =\,a$ & $a>0$ \\ \\
  $\mathbb{A}_{3,3}\oplus \mathbb{A}_1 $& $C_{1,3}^{1}=a,\quad C_{2,3}^{2}=a,\quad C_{3,4}^{1}\,=\,b$& $a>0$ \\
  \hline
\end{tabular}
\end{center}
\begin{center}
table $3$
\end{center}
\end{lemma}
 In this
subsection $(e_1,e_2,e_3,e_4)$ is an orthonormal basis of the Lie
algebra in  Lemma \ref{lemklepikov}.  The matrix of the curvature
operator are given
 in the basis $(e_1 \wedge e_2,e_1 \wedge e_3,e_1
\wedge e_4,e_2 \wedge e_3, e_2 \wedge e_4, e_3 \wedge e_4)$ of
$\wedge^2\mathfrak{g}.$
\begin{enumerate}
\item \textbf{The $Z-$decomposition of the Lie algebra $\mathbb{A}_2\oplus
2\mathbb{A}_1$.}\label{example1}

 For this Lie algebra,
the nonnull components of the curvature tensors are:
\[\begin{array}{lll}
    R(e_1,e_2)e_1=-(a^2+\frac{3}{4}b^2)e_2, &  R(e_1,e_4)e_1=\frac{1}{4}b^2e_4, &  R(e_2,e_4)e_2=\frac{1}{4}b^2e_4, \\
    R(e_1,e_2)e_2=(a^2+\frac{3}{4}b^2)e_1, &  R(e_1,e_4)e_4=-\frac{1}{4}b^2e_1, &  R(e_2,e_4)e_4=-\frac{1}{4}b^2e_2
  \end{array}
.\] and the matrix of the linear curvature operator is
$\texttt{diag}\left(-(a^2+\frac{3}{4}b^2),0,\frac{1}{4}b^2,0,\frac{1}{4}b^2,0\right)$.

\underline{Case I} : $b\neq 0$

 The eigenvectors $h_1=e_1\wedge e_2,\quad h_2=e_1\wedge
e_4,\quad h_3=e_2\wedge e_4$ associated respectively to nonnull
eigenvalues $-(a^2+\frac{3}{4}b^2),\quad \frac{1}{4}b^2,\quad
\frac{1}{4}b^2$\quad are irreducible.
$H_{1}^{0}=\text{Span}\{e_3,e_4\},$\quad
$H_{1}^{1}=\text{Span}\{e_1,e_2\}$,\quad
$H_{2}^{0}=\text{Span}\{e_2,e_3\}$,\quad
$H_{2}^{1}=\text{Span}\{e_1,e_4\}$, \quad
$H_{3}^{0}=\text{Span}\{e_1,e_3\}$
 and $H_{3}^{1}=\text{Span}\{e_2,e_4\}$.

$H_{1}^{1}\nsubseteq H_{2}^{0}$ and $H_{1}^{1}\nsubseteq H_{3}^{0}$.
The  subspace
$V_1=H_{1}^{1}+H_{2}^{1}+H_{3}^{1}=\texttt{Span}\{e_1,e_2,e_4\}$ is
one of subspaces $V_j, j>0$ and the $V-$decomposition of
$\mathbb{A}_2\oplus 2\mathbb{A}_1$\,is\,\, $\mathbb{A}_2\oplus
2\mathbb{A}_1=V_0\oplus V_1,$ where $V_0=\texttt{Span}\{e_3\}$. The
$Z-$decomposition of $\mathbb{A}_2\oplus 2\mathbb{A}_1$ with $b\neq
0$ is
\[\mathbb{A}_2\oplus 2\mathbb{A}_1=Z_0\oplus Z_1
,\]
 where $Z_1=\texttt{Span}\{e_1,e_2,e_4\}$ and
 $Z_0=\texttt{Span}\{e_3\}$

For the Lie algebra $Z_1,$ the Lie brackets are
\[ [e_2,e_1]=-ae_1-be_4, \quad [e_2,e_4]=0,\quad [e_1,e_4]=0 .\]
Therefore, the left invariant Riemannian metric induced by the
restriction on  $Z_1$ of the  inner product $\langle\, , \,\rangle$
is not locally symmetric, see Lemma \ref{lemma2}. In other words,
the metric $g$ is not locally symmetric . Moreover
$[Z_0,Z_1]=\{0\}$.

\underline{Case II} : $b= 0$.

 The eigenvector $h_1=e_1\wedge e_2$ associated  to the  nonnull
eigenvalue $-a^2$ is irreducible. $\quad
H_{1}^{0}=\texttt{Span}\{e_3,e_4\}\quad \texttt{and}\quad
H_{1}^{1}=\texttt{Span}\{e_1,e_2\}$. The subspace
$V_1=\texttt{Span}\{e_1,e_2\}$ is one of the subspaces $V_j,\, j>0$.
The $V-$decomposition is $\mathbb{A}_2\oplus 2\mathbb{A}_1=V_0\oplus
V_1,$ where $V_0=\texttt{Span}\{e_3,e_4\}$.
 The $Z-$decomposition of
$\mathbb{A}_2\oplus 2\mathbb{A}_1$ with $b=0$ is:
\[\mathbb{A}_2\oplus 2\mathbb{A}_1=Z_0\oplus Z_1,\] where
$
 Z_1=\texttt{Span}\{e_1,e_2\}$ and
$Z_0=\texttt{Span}\{e_3,e_4\}$.

The left invariant Riemannian metrics  induced by the restriction on
$Z_i, i\in \{0,1\}$ of the inner product $\langle\, , \,\rangle$ are
locally symmetric,
 see Remark \ref{metricdim2}. Therefore, the metric $g$ is locally symmetric.
Furthermore $Z_0$ is a $2-$dimensional abelian Lie algebra and
$[Z_0,Z_1]=\{0\}$.

\item \textbf{The $Z-$decomposition of the Lie algebra  $\mathbb{A}_{3,1}\oplus \mathbb{A}_1$.}

  The matrix of the  curvature
operator is
$\texttt{diag}(\frac{1}{4}a^2,\frac{1}{4}a^2,0,-\frac{3}{4}a^2
,0,0)$.
 The eigenvectors $h_1=e_1\wedge e_2,\quad h_2=e_1\wedge
e_3,\quad h_3=e_2\wedge e_3$ associated respectively to nonnull
eigenvalues $\frac{1}{4}a^2,\quad \frac{1}{4}a^2,\quad
-\frac{3}{4}a^2$\quad are irreducible. By direct computation,
 the $Z-$decomposition of
$\mathbb{A}_{3,1}\oplus \mathbb{A}_1$  is:
\[\mathbb{A}_{3,1}\oplus \mathbb{A}_1=Z_0\oplus Z_1,\] where
$Z_1=\texttt{Span}\{e_1,e_2,e_3\}$ and
  $Z_0=\texttt{Span}\{e_4\}.$

  For the Lie algebra $Z_1,$ the Lie brackets are
\[ [e_1,e_2]=0, \quad [e_2,e_3]=a e_1,\quad [e_3,e_1]=0. \]
Therefore, by Lemma \ref{lemma1}, the left invariant Riemannian
metric induced by the restriction on $Z_1$ of the inner product
$\langle\, , \,\rangle$ is not locally symmetric.  It follow that
the metric $g$ is not locally symmetric. In addition
$[Z_0,Z_1]=\{0\}$.

\item \textbf{The $Z-$decomposition of the Lie algebra  $\mathbb{A}_{3,3}\oplus
\mathbb{A}_1$.}\label{exemple3}

 The matrix of the linear curvature operator is:
\[\left(
    \begin{array}{cccccc}
      -a^2 & 0 & 0 & 0 &  \frac{-1}{2}ab& 0 \\
      0 & -a^2+\frac{1}{4}b^2 & 0 & 0 & 0 & -ab \\
      0 & 0 & \frac{1}{4}b^2 & 0& 0 & 0 \\
      0 & 0 & 0 & -a^2 & 0 & 0 \\
\frac{-1}{2}ab & 0 & 0 & 0 & 0 & 0 \\
      0 & -ab & 0 & 0 & 0 & -\frac{3}{4}b^2 \\
    \end{array}
  \right)
\]

\underline{Case I}: $b=0$.

The matrix of the linear curvature operator is
$\texttt{diag}(-a^2,\,-a^2,\,0,\,-a^2,\,0,\,0,)$.
 The eigenvectors $h_1=e_1\wedge e_2,\quad h_2=e_1\wedge
e_3,\quad h_3=e_2\wedge e_3$ associated respectively to nonnull
eigenvalues $-a^2,\quad -a^2,\quad -a^2$\quad are irreducible. By
direct computation, the $Z-$decomposition of $\mathbb{A}_{3,3}\oplus
\mathbb{A}_1$  is:
\[\mathbb{A}_{3,3}\oplus \mathbb{A}_1=Z_0\oplus Z_1,\] where
$Z_1=\texttt{Span}\{e_1,e_2,e_3\}$
 $Z_0=\texttt{Span}\{e_4\}$.

 For the Lie algebra $Z_1,$ the Lie brackets are
\[ [e_3,e_1]=-a e_1+0 e_2 , \quad [e_3,e_2]=0 e_1-a e_2,\quad
\text{and} \quad [e_1,e_2]=0. \]
 Therefore, by Lemma \ref{lemma2}, the left invariant Riemannian metric
induced by the restriction on $Z_1$ of the inner product $\langle\,
, \,\rangle$  is locally symmetric. Hence
  $g$ is a locally symmetric. Moreover, $[Z_0,Z_1]=\{0\}$.

\underline{Case II}: $b\neq 0$.

 The eigenvalues , their multiplicities  and the
associated eigenvectors in the basis $(e_1 \wedge e_2,e_1 \wedge
e_3,e_1 \wedge e_4,e_2 \wedge e_3, e_2 \wedge e_4, e_3 \wedge e_4)$
of $\wedge^2\mathfrak{g}$ of the above matrix are contained in the
following table:
\begin{center}
\begin{tabular}{ccl}
  \hline
  Eigenvalue &multiplicity & Eigenvectors \\ \hline
  $1/2\, \left( -a+\sqrt {{a}^{2}+{b}^{2}} \right) a$ & $1 $& $h_1=\dfrac{a-\sqrt {{a}^{2}+{b}^{2}}}{b}e_1 \wedge e_2+e_2 \wedge e_4$
  \\ \\
 $ 1/2\, \left( -a-\sqrt {{a}^{2}+{b}^{2}} \right)a$ & $1$ & $h_2=\dfrac{a+\sqrt {{a}^{2}+{b}^{2}}}{b}e_1 \wedge e_2+e_2 \wedge e_4$
 \\ \\
  $\frac{1}{4}b^2 $& 2 & $h_3=e_1 \wedge
e_4,\,\,h_4=-\dfrac{b}{a}e_1 \wedge e_3+e_3 \wedge e_4
  $\\ \\
  $-{a}^{2}-3/4\,{b
}^{2}$ & $1$ & $h_5=\dfrac{a}{b}e_1 \wedge e_3+e_3 \wedge e_4
  $ \\ \\
 $ -a^2$ & $1$ &$ h_6=e_2 \wedge e_3$ \\
  \hline
\end{tabular}
\end{center}
\begin{center}
Table $4$
\end{center}
By direct computation we prove that\, $\dim\ker\widetilde{h}_i=2$,\,
$i \in \{1,2,3,4,5,6\}$ therefore, they are irreducible  by Remark
\ref{remireducible}, and we also prove by direct computation that,
if $b\neq 0$,\, $\mathbb{A}_{3,3}\oplus \mathbb{A}_1$ does not admit
a $Z-$decomposition.

 \end{enumerate}

\begin{remark}
Let $(G,g)$ be a connected Riemannian Lie group with the associated
euclidian Lie algebra $(\mathbb{A}_{3,3}\oplus \mathbb{A}_1,
\langle\,,\,\rangle)$.
\begin{enumerate}
\item If $b=0$, then the metric $g$ is  locally symmetric.
\item If $b\neq 0$, we can't conclude.
\end{enumerate}
\end{remark}

\subsection{Proof of Theorem \ref{mainthm2}}
O.P.Gladunova  in \cite{gladunodivweyl} classified $4$-dimensional
unimodular harmonic Weyl tensor Riemannian  Lie groups. For
$4$-dimensional nonunimodular Riemannian Lie groups, we refer to
D.S. Voronov  in
 \cite{dsvoronovdivweyl}.  Table $5$ below indicates the Lie algebras and the structure constants in an orthonormal basis
 $(e_1,e_2,e_3,e_4)$
for harmonic Weyl tensor  with $W\neq 0$, for $4-$dimensional
Riemannian Lie groups and the restrictions on their structure
constants.
\begin{center}
Table 5: Lie algebras of $4$-dimensional  $\mathcal{C}-$spaces
Riemannian Lie groups with $W\neq 0$.
\end{center}

 \begin{center}
\begin{tabular}{cll}
  \hline
  Lie algebras & Structure constants & Constraints on $C_{i,j}^{k}$ \\ \hline
  &Decomposable nonunimodular& \\ \\
  $\mathbb{A}_2\oplus 2\mathbb{A}_1$&$C_{1,2}^{2}=a$& $a>0$ \\
  $2\mathbb{A}_{2}$&$C_{1,2}^{2}=a, \,\,C_{3,4}^{4}=b$& $a>0,\,b>0$ \\
   \\
    &Indecomposable nonunimodular& \\ \\
    $\mathbb{A}_{4,5}^{\alpha , \beta}$&$C_{1,4}^{1}=C_{2,4}^{2}=C_{3,4}^{3}=a$&$a > 0$\\
    $\mathbb{A}_{4,6}^{\alpha , \beta}$&$C_{1,4}^{1}=\alpha a,\,\,C_{2,4}^{3}=-C_{3,4}^{2}=-a$&$ \alpha \neq 1,\,\,a > 0$\\
    $\mathbb{A}_{4,6}^{\alpha , \beta}$&$C_{1,4}^{1}=C_{2,4}^{2}=C_{3,4}^{3}=\beta a,\,\,
    C_{2,4}^{3}=-C_{3,4}^{2}=-a$&$ \beta>0,\,\,a > 0$\\
 $\mathbb{A}_{4,9}^{\beta}$&$C_{1,4}^{1}=C_{2,3}^{1}=2a,\,\,C_{2,4}^{2}=C_{3,4}^{3}=a$&$a >
 0$\\ \\
  $\mathbb{A}_{4,11}^{\alpha}$ &$\begin{array}{c}
                                   C_{1,4}^{1}=C_{2,3}^{1}=2a\alpha, \\
                                    C_{2,4}^{2}=C_{3,4}^{3}=a\alpha,\\
C_{2,4}^{3}=-C_{3,4}^{2}=-a
                                 \end{array}
  $ & $\alpha > 0,\,\,a>0$\\
 \hline
\end{tabular}
\end{center}

\vspace{0.5cm}

 We recall the complete list of $4-$dimensional
metric Lie algebras satisfying $W=0$ in  Table $6$, given by E.D.
Rodionov in
    \cite{radionovhalfcon}.  We refer to Proposition
  \ref{gladutheodiago} for the orthonomal
  basis and the spectrum of the curvature operator.

\begin{prof}
According to the classification table of $4-$dimensional Riemannian
Lie groups with nonnull harmonic Weyl tensor, see table $5$, we
compute the matrix of the curvature tensor in the basis $(e_1 \wedge
e_2,e_1 \wedge e_3,e_1 \wedge e_4,e_2 \wedge e_3, e_2 \wedge e_4,
e_3 \wedge e_4)$, \,where $(e_1,e_2,e_3,e_4)$ is an orthonormal
basis of $\mathfrak{g}$.
\begin{enumerate}
\item \textbf{The $Z-$decomposition of the Lie algebra $\mathbb{A}_2\oplus 2\mathbb{A}_1$.}

 See item  (\ref{example1}) of the proof of Theorem
\ref{mainthm1}.

\item \textbf{The $Z-$decomposition of the Lie algebra $2\mathbb{A}_{2}$.}

 The matrix of the curvature operator is
  $ [\overline{R}]=\text{diag}(-a^2,0,0,0,0,-b^2).$
 For the nonnull eigenvalue   $-a^2$ and $-b^2$ of
$[\overline{R}]$, the associated and  irreducible eigenvectors are
$h_1= e_1 \wedge e_2$ and $h_2=e_3 \wedge e_4$  . Therefore
 $H_{1}^{0}=\text{Span}\{e_3,e_4\}$,\quad
 $H_{1}^{1}=\text{Span}\{e_1,e_2\}$,\quad
 $H_{2}^{0}=\text{Span}\{e_1,e_2\}$,\quad \text{and} \quad $H_{2}^{1}=\text{Span}\{e_3,e_3\}$.
\begin{itemize}
 \item[i)] $H_{1}^{1}\subset H_{2}^{0} $, therefore, $H_{1}^{1}$ is one of the subspaces $V_j,\,\,j>0$ in the $V-$decomposition.
\item[ii)] $H_{2}^{1}\subset H_{1}^{0} $, hence, $H_{2}^{1}$ is one of the subspaces $V_i,\,\,i>0$ in the $V-$decomposition.
\end{itemize}
 Setting $V_1=H_{1}^{1}$ and $V_2=H_{2}^{1}$, the $V-$decomposition
of $2\mathbb{A}_{2}$ is $2\mathbb{A}_{2}= V_1\oplus V_2$ and the
$Z-$decomposition is:
\[2 \mathbb{A}_2=Z_1\oplus Z_2,\]
where $Z_1=\texttt{Span}\{e_1,e_2\}$ and
   $Z_2=\texttt{Span}\{e_3,e_4\}.$

The components $Z_1$ and $Z_2$ of the $Z-$decomposition satisfy
$[Z_1,Z_2]=\{0\}$. Moreover, the left invariant Riemannian metrics
induced by the restriction on on $Z_i, i\in \{1,2\}$ of the inner
product $\langle\, , \,\rangle$ are locally symmetric,
 by Remark \ref{metricdim2}.
Thus, the metric $g$ is locally symmetric.

\item \textbf{The $Z-$decomposition of the Lie algebra $\mathbb{A}_{4,6}^{\alpha,\beta},\,\,\alpha \neq
1 $.}

 The nonnull structure constants are $C_{1,4}^{1}=\alpha a, \,\,C_{2,4}^{3}=-a,$ \quad $C_{3,4}^{2}=a,\,a>0$
and the matrix of the curvature tensor linear operator is
   $[\overline{R}]=\texttt{diag}(0,0,-(\alpha a)^2,0,0,0).$
 For the nonnull eigenvalue $-(\alpha a)^2$ of
 $[\overline{R}]$,\quad
 the associated irreducible eigenvector is $h_1= e_1 \wedge e_4$.
  The $Z-$decomposition is
\[\mathfrak{g}=Z_0\oplus Z_1,\]
 where $ Z_1=\text{Span}\{e_1,e_4\}$ and
$Z_0=\text{Span}\{e_2,e_3\}$. On the other hand, the left invariant
Riemannian metrics induced by the restriction on $Z_i, i\in \{0,1\}$
of the inner product $\langle\, , \,\rangle$  are locally symmetric,
 by Remark \ref{metricdim2}. Therefore, the metric $g$ is locally symmetric.  Also $[Z_0,Z_1] = Z_0$ and
$Z_0$ is an abelian Lie algebra.

\item \textbf{The $Z-$decomposition of the Lie algebras  $\mathbb{A}_{3,9}\oplus \mathbb{A}_1$,\quad
$\mathbb{A}_{3,3}\oplus \mathbb{A}_1$,\quad,
 $\mathbb{A}_{3,7}^{\alpha}\oplus \mathbb{A}_1$, \quad \text{and}\quad $\mathbb{A}_{4,12}.$}

  For these Lie algebras, the irreducible eigenvectors associated
  to nonnull eigenvalues are $h_1=e_1 \wedge e_2,\quad h_2=e_1 \wedge
  e_3$ and $h_3=e_2 \wedge e_3.$ By direct computation, the
  $Z-$decomposition is
   \[\mathfrak{g}=Z_0\oplus Z_1,\]
   where $Z_0=\{e_4\}$ and $Z_1=\text{Span}\{e_1,e_2,e_3\}$.
\end{enumerate}

For more details  see the
   case $\mathbb{A}_{3,3}\oplus \mathbb{A}_1$
  in  the proof of Theorem
\ref{mainthm1}.
\end{prof}

\subsubsection{The structure of the $3$-dimensional component $Z_1$.}

The Lie brackets on the $3-$dimensional components $Z_1$ of the
above Lie algebras are:
\begin{enumerate}
\item For $\mathbb{A}_{3,9}\oplus \mathbb{A}_1$,
\[ [e_1,e_2]=-a \sqrt{1+m^2}\, e_3,\quad \,[e_2,e_3]=-a \sqrt{1+m^2}\, e_1,\,\quad [e_3,e_1]=-a \sqrt{1+m^2}\, e_2. \]
   Moreover, $[Z_0,Z_1] \subset Z_1$.

\item For $\mathbb{A}_{3,7}^{\alpha}\oplus
\mathbb{A}_1$,
\[ \,[e_3,e_1]=-\alpha a\, e_1 +a\,e_2,\quad \,[e_3,e_2]=-a  e_1 - \alpha \,a e_2,\,\quad [e_1,e_2]=0. \]
 and it holds that $[Z_0,Z_1] = \{0\}$.

\item For $\mathbb{A}_{3,3}\oplus \mathbb{A}_1$, see
item (\ref{exemple3})  of the proof of Theorem \ref{mainthm1}.
\item For  $\mathbb{A}_{4,12}$,
\[ \,[e_3,e_1]=-\sqrt{a^2+b^2}\, e_1 +\frac{bd}{\sqrt{a^2+b^2}}\,e_2,\quad \,
[e_3,e_2]=-\frac{bd}{\sqrt{a^2+b^2}}\, e_1  -\sqrt{a^2+b^2} e_2,\]
\[ [e_1,e_2]=0. \]
 We also have $[Z_0,Z_1] \subset Z_1$.

For each of these Lie algebras, the left invariant Riemannian metric
induced by the restriction of the inner product $\langle\,
,\,\rangle$ on $Z_1$ are locally symmetric by Lemmas \ref{lemma1}
and \ref{lemma2}. Therefore, the metrics $g$ are locally symmetric.
\end{enumerate}
\begin{remark}
The other Lie algebras in tables $5 \ and \ 6$ do not admit
$Z-$decomposition.
\begin{enumerate}
\item  For the Lie algebras $4\mathbb{A}_1,\quad
\mathbb{A}_{3,6}\oplus \mathbb{A}_1$,\,\,$0$ is the only eigenvalue,
see Remark \ref{vdecopremark}.

\item For the Lie algebras $\mathbb{A}_{4,9}^{\beta}\,,\beta=1,\quad
\mathbb{A}_{4,11}^{\alpha,},\alpha >0$ , $\ker
\widetilde{h}_i=\{0\}$, where $h_i$ is an irreducible eigenvector
associated to nonnull eigenvalue of curvature operator
$\overline{R}$, see Remark \ref{vdecopremark}. For these Lie
algebras, the matrix $\overline{R}$ in the basis $(e_i \wedge
e_j)_{i<j}$ is given by O.P. Gladunova  in \cite{gladunovarusse}.

\item For the Lie algebras $\mathbb{A}_{4,5}^{\alpha,\beta},\quad
\mathbb{A}_{4,6}^{\alpha,\beta},\alpha=\beta$ with $
C_{1,4}^{1}=C_{2,4}^{2}=C_{3,4}^{3}=\alpha
 a, \,\,\text{and}\,\,\,
C_{3,4}^{2}=-C_{2,4}^{3}=a $, by the direct application of
$Z-$decomposition technic, we have $Z_1=\mathfrak{g}$.
\end{enumerate}
\end{remark}
\begin{center}
 Table 6: Lie algebras of $4$-dimensional conformally fat Riemannian Lie groups.
 \end{center}
  \begin{center}
  \begin{tabular}{lll}
    \hline
     Lie algebras & Structure constants & $\text{Spec}(\overline{R})=\left\{K_{12},K_{13},K_{14},K_{23},K_{24},K_{34}\right\}$ \\ \hline
     $4\mathbb{A}_1$ &$C_{i,j}^{k}=0,\,\, \forall i,j,k,$  &$\{0,0,0,0,0,0\}$
     \\\\
     $\begin{array}{l}
             \mathbb{A}_{3,6}\oplus \mathbb{A}_1 \\
             c>0
           \end{array}
     $ & $C_{2,3}^{1}=c,\,\,\,C_{1,3}^{2}=-c$ & $\{0,0,0,0,0,0\}$
     \\\\
     $\begin{array}{l}
            \mathbb{A}_{3,9}\oplus \mathbb{A}_1 \\
            a>0
          \end{array}
    $ &$\begin{array}{l}
                                                   C_{1,3}^{2}=a\sqrt{1+m^2}\\
                                                   C_{1,2}^{3}=C_{2,3}^{1}=-a\sqrt{1+m^2}\\
                                                   C_{2,4}^{1}=-C_{1,4}^{2}=a\,m\sqrt{1+m^2}\\
                                                 \end{array}
    $  & $\left\{\frac{a^2(1+m^2)}{4},\frac{a^2(1+m^2)}{4},0,\frac{a^2(1+m^2)}{4},0,0\right\}$ \\\\
      $
     \begin{array}{l}
        \mathbb{A}_{3,3}\oplus \mathbb{A}_1\\
       a>0
     \end{array}
     $ & $C_{1,3}^{1}=C_{2,3}^{2}=a$ & $\{-a^2,-a^2,0,-a^2,0,0\}$
     \\\\
      $
     \begin{array}{l}
       \mathbb{A}_{3,7}^{\alpha}\oplus \mathbb{A}_1 \\
a>0, \,\,\alpha >0
     \end{array}
     $ &  $
     \begin{array}{l}
       C_{1,3}^{1}=C_{2,3}^{2}=\alpha a \\
       C_{2,3}^{1}=-C_{1,3}^{2}= a
     \end{array}
     $& $\{-\alpha^2\,a^2,-\alpha^2\,a^2,0,-\alpha^2\,a^2,0,0\}$ \\\\
      $
     \begin{array}{l}
       \mathbb{A}_{4,5}^{\alpha,\beta} \\
a > 0\\
\alpha=\beta=1
     \end{array}
     $ &$C_{1,4}^{1}=C_{2,4}^{2}=C_{3,4}^{3}=a$  &
      $\{-a^2,-a^2,-a^2,-a^2,-a^2,-a^2\}$ \\\\
      $\begin{array}{l}
             \mathbb{A}_{4,6}^{\alpha,\beta} \\
\alpha = \beta,\\
\,a>0,\alpha >0
           \end{array}
     $ &$\begin{array}{l}
                                                 C_{1,4}^{1}=C_{2,4}^{2}=C_{3,4}^{3}=\alpha a \\
                                                 C_{3,4}^{2}=-C_{2,4}^{3}=a
                                               \end{array}
     $  & $\{-\alpha^2a^2,-\alpha^2a^2,-\alpha^2a^2,-\alpha^2a^2,-\alpha^2a^2,-\alpha^2a^2\}$ \\\\
     $\begin{array}{l}
             \mathbb{A}_{4,12} \\
a>0,d>0
           \end{array}
    $ & $\begin{array}{l}
                                    C_{1,3}^{1}=C_{2,3}^{2}=\sqrt{a^2+b^2}\\
                                   C_{2,4}^{1}=-C_{1,4}^{2}=\frac{a\,d}{\sqrt{a^2+b^2}} \\
                                   C_{2,3}^{1}=-C_{1,3}^{2}=\frac{b\,d}{\sqrt{a^2+b^2}}
                                 \end{array}
    $ & $\{-(a^2+b^2),-(a^2+b^2),0,-(a^2+b^2),0,0\}$ \\
    \hline
  \end{tabular}
  \end{center}


\begin{thebibliography}{99}
\bibitem{antonioscala} Antonio J. DI scala, Autoparallel
Distribution and Splitting Theorems, \emph{TSUKUBA J. Maths}, \bf{27
(1)} (2003), 99-102.

\bibitem{albesse} Arthur L. Besse, \emph{Einstein Manifold},
Springer-Verlag Berlin Heildelberg (Berlin,1987).

\bibitem{dsvoronovdivweyl} D. S. Voronov and E.D Rodionov, Left-Invariant Riemannian Metrics on
Four-Dimensional Nonunimodular Lie Groups with Zero-Divergence Weyl
Tensor, \emph{Doklady Mathematics}, \bf{81(3)} (2010), 392-394.
\bibitem{Ecartan1} E. Cartan, Sur les classes remarquables d'espaces
Riemanniens, \emph{Bulletin de la S.M.F}, \bf{54} (1926),  214-264.

\bibitem{Ecartan2} E. Cartan, Sur les classes remarquables d'espaces
Riemanniens II, \emph{Bulletin de la S.M.F}, \bf{55} (1927),
114-134.

\bibitem{radionovhalfcon}  E. D. Rodionov, V.V. Slavskii, O. P. Khromova,
On the Curvature Operator Spectrum of (Half)Conformally Flat
Riemannian Metrics, \emph{The news of Altai State University},(2015)
DOI: 10.14258/izvasu(2015)1.1-19.

\bibitem{Jmilnor} John Milnor, Curvature of Left Invariant Metrics on
Lie  Groups, \emph{Adv.Math.} \bf{21}(1976), 293-329.

\bibitem{gladunoconf} O. P. Gladunova, E. D. Rodianov, and V.V.
Slavskii, On the Spectrum of Curvature Operator of Conformally Flat
Riemannian Manifold, \emph{Doklady Mathematics}, \bf{87(3)} (2013),
279-281.

\bibitem{gladunovarusse} O. P. Gladunova, D.N. Oskorbin, An Application of Symbolic Computation
Packages to the Investigation
 of the Curvature Operator Spectrum on the Metric Lie Groups, \emph{The news of Altai State University}, \bf{1}(2013), 19-23.

\bibitem{gladunodivweyl} O. P. Gladunova,  V.V. Slavskii, Harmonicity of the Weyl Tensor of Left-Invariant
Riemannian Metric on Four-Dimensional Unimodular Lie Groups,
\emph{Doklady Mathematics}, \bf{81(2)} (2010), 298-300.

\bibitem{klepikov}  P. N.  Klepikov, D.N. Oskorbin, Milnor's Generalized Bases for Some Four-dimensional Metric Lie
Algebras, \emph{The news of Altai State University},(2015) DOI
10.14258/izvasu(2015)1.1-13.

\bibitem{kobayashi} Shoshichi Kobayashi and Katsumi Nomizu,
\emph{Foundations of Differential Geometry,} Vol I, Wiley, New York,
(1963).

\bibitem{szabolocal} Z. I. Szab\'o, Structure theorems on Riemannian
spaces satisfying $R(X,Y).R\,\,=\,\,0$. The local version, \emph{J.
Differential Geometry}  \bf{17} (1982), 531-582.


\end{thebibliography}
\end{document}